\documentclass[graybox]{svmult}

\usepackage[english]{babel}
\usepackage{mathptmx}
\usepackage{helvet}
\usepackage{courier}
\usepackage{type1cm}
\usepackage{makeidx}
\usepackage{graphicx}
\usepackage{multicol}
\usepackage[bottom]{footmisc}
\usepackage{amsfonts}
\usepackage{amsmath}

\usepackage{hyperref}
\usepackage{pstricks}
\usepackage{pst-circ}

\makeindex

\newcommand{\diag}{\mathop{\text{diag}}}
\newcommand{\divergence}{\mathop{\text{div}}}
\newcommand{\spanvec}{\mathop{\text{span}}}
\DeclareMathOperator*{\argmax}{arg\,max}
\DeclareMathOperator*{\argmin}{arg\,min}
\newcommand{\R}{\mathbb{R}}
\newcommand{\N}{\mathbb{N}}
\newcommand{\RT}{\text{RT}}
\newcommand{\calP}{\mathcal{P}}
\newtheorem{algorithm}[theorem]{Algorithm}

\begin{document}

\title*{Residual Based Sampling in POD Model Order Reduction of Drift-Diffusion Equations in Parametrized Electrical Networks}
\titlerunning{Residual Based Sampling in POD Model Order Reduction}
\author{Michael Hinze and Martin Kunkel}
\institute{
     Michael Hinze \and Martin Kunkel \at
     Department of Mathematics, University of Hamburg, Bundesstr. 55, 20146 Hamburg, Germany\\ email: \url{michael.hinze@uni-hamburg.de},\ \ \url{martin.kunkel@uni-hamburg.de}
}

\maketitle

\abstract{We consider integrated circuits with semiconductors modeled by modified nodal analysis and drift-diffusion equations. The drift-diffusion equations are discretized in space using mixed finite element method. This discretization yields a high dimensional differential-algebraic equation. We show how proper orthogonal decomposition (POD) can be used to reduce the dimension of the model.
We compare reduced and fine models and give numerical results for a basic network with one diode. Furthermore we discuss an adaptive approach to construct POD models which are valid over certain parameter ranges. Finally, numerical investigations for the reduction of a 4-diode rectifier network are presented, which clearly indicate that POD model reduction delivers surrogate models for the diodes involved, which depend on the position of the semiconductor in the network.}

\bigskip


\section{Introduction} \label{sec:introduction}
\begin{figure} \sidecaption
\includegraphics[width=0.55\textwidth]{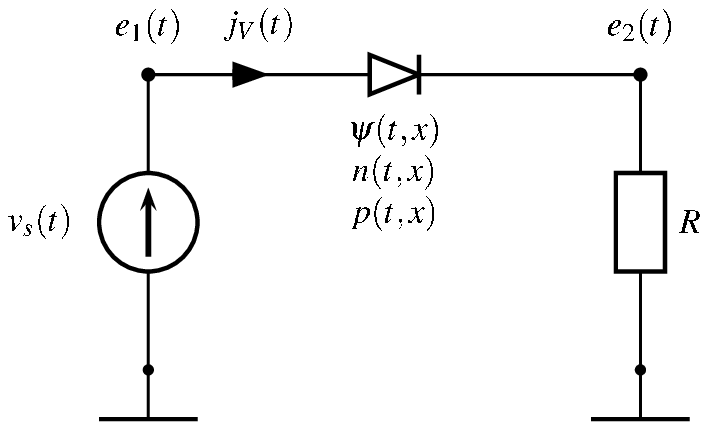}
\caption[Basic test circuit with one diode.]{Basic test circuit with one diode. The network is described by
\begin{align*}
    A_V  &= \begin{pmatrix}
                    \phantom{-}1, & \phantom{-}0
              \end{pmatrix}^\top, \\
    A_S  &= \begin{pmatrix}
                     -1, &  \phantom{-}1
             \end{pmatrix}^\top, \\
    A_R  &= \begin{pmatrix}
                    \phantom{-}0, & \phantom{-}1
              \end{pmatrix}^\top, \\
    g(A_R^\top e,t) &= \frac{1}{R} e_2(t).
\end{align*}}
\label{fig:basic}
\end{figure}

In this article we investigate a POD-based model order reduction for semiconductors in electrical networks. Electrical networks can be efficiently modeled by a differential-algebraic equation (DAE) which is obtained from modified nodal analysis. Denoting by $e$ the node potentials and by $j_L$ and $j_V$ the currents of inductive and voltage source branches, the DAE reads (see \cite{GFM05,Tis03})
    \begin{align}
        A_C \frac{d}{dt} q_C(A_C^\top e, t) + A_R g(A_R^\top e, t) + A_L j_L + A_V j_V &= - A_I i_s(t),
          \label{eq:mna:1} \\
        \frac{d}{dt} \phi_L(j_L, t) - A_L^\top e &= 0 \label{eq:mna:2}, \\
        A_V^\top e &= v_s(t). \label{eq:mna:3}
    \end{align}
Here, the incidence matrix $A = [A_R,A_C,A_L,A_V,A_I]$ represents the network topology, e.g. at each non mass node $i$, $a_{ij} = 1$ if the branch $j$ leaves node $i$ and $a_{ij} = -1$ if the branch $j$ enters node $i$. If node $i$ and branch $j$ are not connected, then $a_{ij} = 0$. $q_C$, $g$ and $\phi_L$ are continuously differentiable functions defining the voltage-current relations of the network components. The continuous functions $v_s$ and $i_s$ are the voltage and current sources. For a basic example consider the network in Figure~\ref{fig:basic}. Under the assumption that the Jacobians
    \[ D_C(e,t) := \frac{\partial q_C}{\partial e}(e,t), \quad
       D_G(e,t) := \frac{\partial g}{\partial e}(e,t), \quad
       D_L(j,t) := \frac{\partial \phi_L}{\partial j}(j,t) \]
are positive definite, analytical properties (e.g. the index) of DAE \eqref{eq:mna:1}-\eqref{eq:mna:3} are investigated in \cite{EFM03} and \cite{ET00}. In linear networks, the matrices $D_C$, $D_G$ and $D_L$ are positive definite diagonal matrices with capacitances, conductivities and inductances on the diagonal.

Often semiconductors themselves are modeled by electrical networks. These models are stored in a library and are stamped into the surrounding network in order to create a complete model of the integrated circuit. Here we use a different approach which models semiconductors by the transient drift-diffusion equations. Advantages of this approach are the higher accuracy of the model and fewer model parameters. On the other hand, numerical simulations are more expensive. For a comprehensive overview of the drift-diffusion equations we refer to \cite{AAR03,ANR05,BMP05,Mar86,Sel84}. Using the notation introduced there, we have the following system of partial differential equations for the electrostatic potential $\psi(t,x)$, the electron and hole concentrations $n(t,x)$ and $p(t,x)$, and the current densities $J_n(t,x)$ and $J_p(t,x)$:
    \begin{align*}
        \divergence (\varepsilon \nabla \psi) &= q(n-p-C), \\
        -q \partial_t n + \divergence J_n &= \phantom{-}q R(n,p,J_n,J_p), \\
        \phantom{-} q \partial_t p + \divergence J_p &= -q R(n,p,J_n,J_p),\\
        J_n &= \mu_n q( \phantom{-}U_T \nabla n - n \nabla \psi), \\
        J_p &= \mu_p q( -U_T \nabla p - p \nabla \psi),
    \end{align*}
with $(t,x) \in [0,T] \times \Omega$ and $\Omega \subset \R^d$ ($d=1,2,3$). The nonlinear function $R$ describes the rate of electron/hole recombination, $q$ is the elementary charge, $\varepsilon$ the dielectricity, $\mu_n$ and $\mu_p$ are the mobilities of electrons and holes. The temperature is assumed to be constant which leads to a constant thermal voltage $U_T$. The function $C$ is the time independent doping profile.

This paper is organized as follows. In Section~\ref{sec:fullsystem} we present the (scaled) mathematical model for the complete coupled network including semiconductors modeled by the DD equations. In Section~\ref{sec:fullsim} we describe the numerical treatment of the network using the method of lines. For the spacial discretization of the DD models we use mixed finite elements. For the time integration of the resulting DAE system we use DASSL with step width control and so obtain snapshots $y^h(t_i,\cdot)$, $i=1,\ldots,k$, which represent the state of the circuit and the semiconductors at time $t_i$. Based on these snapshots and POD we construct a reduced model in Section~\ref{sec:red}. We discuss the properties of the reduced model with respect to parameter changes in Section~\ref{sec:param}. In order to obtain a reduced model which is valid over a considerable range of parameters it is necessary to refine the model. For this purpose we adapt the reduced basis method combined with the greedy approach proposed in~\cite{PR07} to our setting. In Section~\ref{sec:numerics} we finally present numerical experiments, and also discuss advantages and shortcomings of our approach.

\section{Complete coupled system} \label{sec:fullsystem}
In the present section we develop the complete coupled system for a network with $n_s$ semiconductors. We will not specify an extra index for semiconductors, but we keep in mind that all semiconductor equations and coupling conditions need to be introduced for each semiconductor.

\begin{figure} \sidecaption
\includegraphics[width=0.6\textwidth]{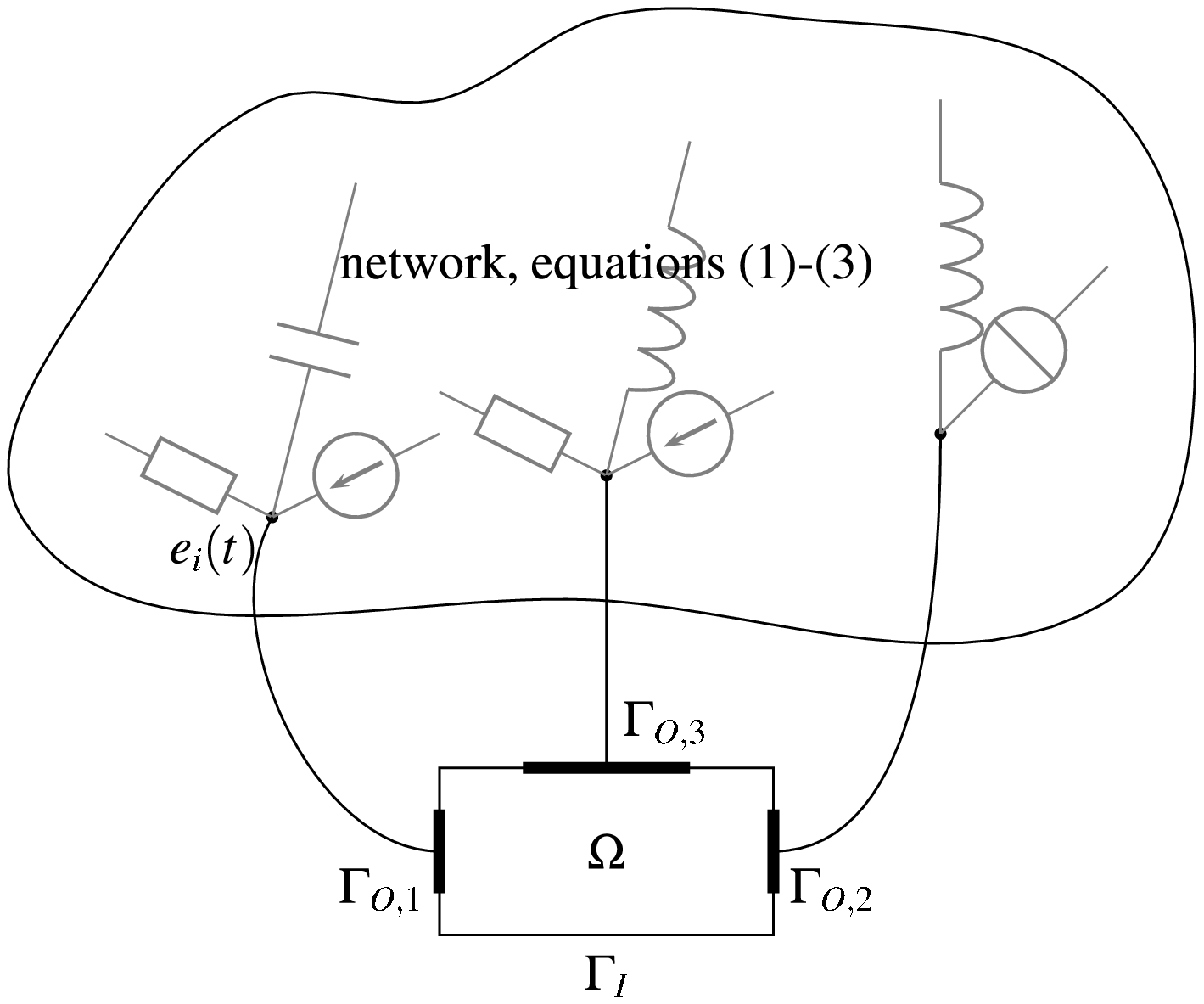}
\caption[Sketch of a coupled system.]{Sketch of a coupled system with one semiconductor. Here
\[ \psi(t,x) = e_i(t) + \psi_{bi}(x), \]
for all $(t,x) \in [0,T] \times \Gamma_{O,1}$.}
\label{fig:sketch}
\end{figure}

For the sake of simplicity we assume that to a semiconductor $m$ semiconductor interfaces $\Gamma_{O,k} \subseteq \Gamma \subset \partial \Omega$,  $k=1,\ldots,m$ are associated, which are all Ohmic contacts, compare Figure~\ref{fig:sketch}. The dielectricity $\varepsilon$ shall be constant over the whole domain $\Omega$. We focus on the Shockley-Read-Hall recombination
  \[ R(n,p) := \frac{n p -\eta^2}{\tau_p(n+\eta) + \tau_n(p+\eta)} \]
which does not depend on the current densities. Herein, $\tau_n$ and $\tau_p$ are the average lifetimes of electrons and holes and $\eta$ is the constant intrinsic concentration which satisfy $\eta^2 = n p$ if the semiconductor is in thermal equilibrium.

The coupling of the drift-diffusion equations and the electrical network is established as follows. The current $A_S j_S$ is added to equation~\eqref{eq:mna:1}, where
\begin{equation}
    j_{S,k} = \int_{\Gamma_{O,k}} (J_n+J_p - \varepsilon \partial_t \nabla \psi) \cdot \nu\, d\sigma.
    \label{eq:current}
\end{equation}
I.e. the current is the integral over the current density $J_n + J_p$ plus the displacement current in normal direction $\nu$. Furthermore, the node potentials of nodes which are connected to a semiconductor interface are introduced in the boundary conditions of the drift-diffusion equations (see also Figure~\ref{fig:sketch}):
    \begin{align}
        \psi(t,x) &= \psi_{bi}(x) + (A_S^\top e(t))_k = U_T \log \left( \frac{\sqrt{ C(x)^2 + 4 \eta^2 } + C(x)}{2 \eta} \right) + (A_S^\top e(t))_k, \label{eq:boundary:psi}\\
        n(t,x)    &= \frac{1}{2} \left( \sqrt{ C(x)^2 + 4 \eta^2 } + C(x) \right), \label{eq:boundary:n} \\
        p(t,x)    &= \frac{1}{2} \left( \sqrt{ C(x)^2 + 4 \eta^2 } - C(x) \right), \label{eq:boundary:p}
    \end{align}
for $(t,x) \in [0,T] \times \Gamma_{O,k}$. Here, $\psi_{bi}(x)$ denotes the build-in potential and $\eta$ the constant intrinsic concentration. All other parts of the boundary are isolation boundaries $\Gamma_I := \Gamma \setminus \Gamma_O$, where $\nabla \psi \cdot \nu = 0$, $J_n \cdot \nu = 0$ and $J_p \cdot \nu = 0$ holds.

The complete model forms a partial differential-algebraic equation (PDAE). The analytical and numerical analysis of such systems is subject to current research, see \cite{BT07, Gue01, ST05, Tis03}. The simulation of the complete coupled system is expensive and numerically difficult due to bad scaling of the drift-diffusion equations. The numerical issues can be significantly reduced by the unit scaling procedure discussed in \cite{Sel84}. That means we substitute
\begin{multline*}
  x = L \tilde{x}, \quad
  \psi = U_T \tilde{\psi}, \quad
  n = \|C\|_\infty \tilde{n}, \quad
  p = \|C\|_\infty \tilde{p}, \quad
  C = \|C\|_\infty \tilde{C}, \\
  J_n = \frac{q U_T \|C\|_\infty}{L} \mu_n \tilde{J_n}, \quad
  J_p = \frac{q U_T \|C\|_\infty}{L} \mu_p \tilde{J_p}, \quad \eta = \tilde{\eta} \|C\|_\infty,
\end{multline*}
where $L$ denotes a specific length of the semiconductor. The scaled drift-diffusion equations then read
    \begin{align}
        \lambda \Delta \psi &= n-p-C, \label{eq:dd:psi} \\
        - \partial_t n + \nu_n \divergence J_n &= \phantom{-}R(n,p), \label{eq:dd:n}\\
          \partial_t p + \nu_p \divergence  J_p &= - R(n,p),\label{eq:dd:p} \\
         J_n &=  \phantom{-} \nabla n - n \nabla \psi, \label{eq:dd:Jn} \\
         J_p &= - \nabla p - p \nabla \psi, \label{eq:dd:Jp}
    \end{align}
where we omit the tilde for the scaled variables. The constants are given by $\lambda := \frac{\varepsilon U_T}{L^2 q \|C\|_\infty}$, $\nu_n := \frac{U_T \mu_n}{L^2}$ and $\nu_p := \frac{U_T \mu_p}{L^2}$.

\section{Simulation of the full system} \label{sec:fullsim}
Classical approaches for the simulation of drift-diffusion equations (e.g. Gummel iterations \cite{Gum64}) approximate $J_n$ and $J_p$ by piecewise constant functions and then solve equations \eqref{eq:dd:n} and \eqref{eq:dd:p} with respect to $n$ and $p$ explicitly. This helps reducing the computational effort and increases the numerical stability. For the model order reduction approach proposed in the present work this method has the disadvantage of introducing additional non-linearities, arising from the exponential structure of the Slotboom variables, see \cite{ST05}.

Since the electrical field represented by the (negative) gradient of the electrical potential $\psi$ plays a dominant role in~\eqref{eq:dd:psi}-\eqref{eq:dd:Jp} and is present also in the coupling condition~\eqref{eq:current}, we provide for it the additional variable $g_\psi = \nabla \psi$ leading to the following mixed formulation of the DD equations:
    \begin{align}
        \lambda \divergence g_\psi &= n-p-C, \label{eq:mixed:psi} \\
        - \partial_t n + \nu_n \divergence J_n &= \phantom{-}R(n,p) \label{eq:mixed:n},\\
          \partial_t p + \nu_p \divergence  J_p &= - R(n,p) \label{eq:mixed:p},\\
         g_\psi &=  \phantom{-} \nabla \psi, \label{eq:mixed:gpsi} \\
         J_n &=  \phantom{-} \nabla n - n g_\psi, \label{eq:mixed:Jn} \\
         J_p &= - \nabla p - p g_\psi. \label{eq:mixed:Jp}
    \end{align}

The weak formulation of \eqref{eq:mixed:psi}-\eqref{eq:mixed:Jp} then reads: Find $\psi, n, p \in [0,T] \times L^2(\Omega)$ and $g_\psi, J_n, J_p \in [0,T] \times H_{0,N}(\divergence, \Omega)$ such that
    \begin{align}
         \lambda \int_\Omega \divergence g_\psi\ \varphi &= \phantom{-}\int_\Omega (n-p)\ \varphi - \int_\Omega C \ \varphi,\label{eq:weak:psi}\\
        -\int_\Omega \partial_t n\ \varphi + \nu_n \int_\Omega \divergence J_n\ \varphi &= \phantom{-}\int_\Omega R(n,p)\ \varphi,\label{eq:weak:n}\\
         \int_\Omega \partial_t p\ \varphi + \nu_p \int_\Omega \divergence J_p\ \varphi &= -\int_\Omega R(n,p)\ \varphi,\label{eq:weak:p}\\
        \int_\Omega  g_\psi\cdot \phi &= -\int_\Omega \psi\ \divergence \phi + \int_\Gamma \psi\ \phi \cdot \nu,\label{eq:weak:gpsi} \\
        \int_\Omega  J_n\cdot \phi &= -\int_\Omega n\ \divergence \phi + \int_\Gamma n\ \phi \cdot \nu - \int_\Omega n\ g_\psi \cdot \phi,\label{eq:weak:Jn} \\
        \int_\Omega  J_p\cdot \phi &= \phantom{-}\int_\Omega p\ \divergence \phi - \int_\Gamma p\ \phi \cdot \nu - \int_\Omega p\ g_\psi \cdot \phi,\label{eq:weak:Jp}
    \end{align}
are satisfied for all $\varphi \in L^2(\Omega)$ and $\phi \in H_{0,N}(\divergence, \Omega)$ where the space $H_{0,N}(\divergence, \Omega)$ is defined by
\begin{align*}
    H(\divergence, \Omega) &:= \{ y \in L^2(\Omega)^d:\ \divergence y \in L^2(\Omega) \},\\
    H_{0,N}(\divergence, \Omega) &:= \left\{ y \in H(\divergence,\Omega):\ y \cdot \nu = 0 \text{ on } \Gamma_I \right\}.
\end{align*}
Consequently, the boundary integrals on the right hand sides in equations \eqref{eq:weak:gpsi}-\eqref{eq:weak:Jp} reduce to integrals over the interfaces $\Gamma_{O,k}$, where the values of $\psi$, $n$ and $p$ are determined by the Dirichlet boundary conditions~\eqref{eq:boundary:psi}-\eqref{eq:boundary:p}. We note that, in contrast to the standard weak form associated with~\eqref{eq:dd:psi}-\eqref{eq:dd:Jp}, the Dirichlet boundary values are naturally included in the weak formulation~\eqref{eq:weak:psi}-\eqref{eq:weak:Jp} and the Neumann boundary conditions have to be included in the space definitions. This is advantageous in the context of POD model order reduction since the non-homogeneous boundary conditions~\eqref{eq:boundary:psi}-\eqref{eq:boundary:p} are not present in the space definitions.

Here, equations~\eqref{eq:weak:psi}-\eqref{eq:weak:Jp} are discretized in space with Raviart-Thomas finite elements of degree $0$ ($RT_0$), alternative discretization schemes for the mixed problem are presented in~\cite{BMP05}. To describe the $RT_0$-approach for $d=2$ spatial dimensions, let $\mathcal{T}$ be a triangulation of $\Omega$ and let $\mathcal{E}$ be the set of all edges. Let $\mathcal{E}_I := \{ E \in \mathcal{E}:\ E \subset \bar\Gamma_I\}$ be the set of edges at the isolation (Neumann) boundaries. The potential and the concentrations are approximated in space by piecewise constant functions
    \[ \psi^h(t), n^h(t), p^h(t) \in L_h := \{ y \in L^2(\Omega):\ y|_T(x) = c_T,\ \forall T \in \mathcal{T} \}, \]
and the discrete fluxes $g^h_\psi(t)$, $J^h_n(t)$ and $J^h_p(t)$ are elements of the space
\begin{multline*}
    RT_0 := \{ y : \Omega \rightarrow \R^d:\ y|_T(x) = a_T + b_T x,\ a_T \in \R^d,\ b_T \in \R,\\
        [y]_E \cdot \nu_E = 0,\ \text{for all inner edges } E \}.
\end{multline*}
Here, $[y]_E$ denotes the jump $y|_{T_+} - y|_{T_-}$ over a shared edge $E$ of the elements $T_+$ and $T_-$. The continuity assumption yields $\RT_0 \subset H(\divergence,\Omega)$.

We set
    \[H_{h,0,N}(\divergence,\Omega) := \left(\RT_0 \cap H_{0,N}(\divergence,\Omega) \right) \subset H_{0,N}(\divergence,\Omega).\]
Then it can be shown, that $H_{h,0,N}$ posses an edge-oriented basis $\{\phi_j\}_{j=1,\ldots,M}$. We use the following finite element Ansatz in~\eqref{eq:weak:psi}-\eqref{eq:weak:Jp}
\begin{equation}\label{eq:galerkinansatz}
    \left.
    \begin{aligned}
      \psi^h(t,x) &= \sum_{i=1}^N \psi_i(t) \varphi_i(x), \quad
      &g_\psi^h(t,x) &= \sum_{j=1}^M g_{\psi,j}(t) \phi_j(x), \\
      n^h(t,x) &= \sum_{i=1}^N n_i(t) \varphi_i(x), \quad
      &J_n^h(t,x) &= \sum_{j=1}^M J_{n,j}(t) \phi_j(x),\\
      p^h(t,x) &= \sum_{i=1}^N p_i(t) \varphi_i(x), \quad
      &J_p^h(t,x) &= \sum_{j=1}^M J_{p,j}(t) \phi_j(x),
    \end{aligned}
    \quad \right\}
\end{equation}
where $N := |\mathcal{T}|$, i.e. the number of elements of $\mathcal{T}$, and $M :=|\mathcal{E}| - |\mathcal{E}_N|$, i.e. the number of inner and Dirichlet boundary edges.

This in~\eqref{eq:weak:psi}-\eqref{eq:weak:Jp} yields
\begin{align*}
         \lambda \sum_{j=1}^M g_{\psi,j}(t) \int_\Omega \divergence \phi_j\ \varphi_k - \sum_{i=1}^N (n_i(t)-p_i(t))  \int_\Omega \varphi_i\ \varphi_k &=-\int_\Omega C \ \varphi_k,\\
        -\sum_{i=1}^N \dot{n}_i(t) \int_\Omega \varphi_i\ \varphi_k + \nu_n \sum_{j=1}^M J_{n,j}(t)\int_\Omega \divergence \phi_j\ \varphi_k -\int_\Omega R(n^h, p^h)\ \varphi_k&=0,\\
         \sum_{i=1}^N \dot{p}_i(t) \int_\Omega \varphi_i\ \varphi_k + \nu_p \sum_{j=1}^M J_{p,j}(t)\int_\Omega \divergence \phi_j\ \varphi_k +\int_\Omega R(n^h, p^h)\ \varphi_k&=0,\\
         \sum_{j=1}^M g_{\psi,j}(t) \int_\Omega \phi_j\cdot \phi_l +\sum_{i=1}^N \psi_i(t) \int_\Omega \varphi_i\ \divergence \phi_l &=\phantom{-} \int_\Gamma \psi^h\ \phi_l \cdot \nu, \\
        \sum_{j=1}^M J_{n,j}(t) \int_\Omega  \phi_j\cdot \phi_l +\sum_{i=1}^N n_i(t) \int_\Omega \varphi_i\ \divergence \phi_l + \int_\Omega n^h  g^h_{\psi} \cdot \phi_l&=\phantom{-}\int_\Gamma n^h\ \phi_l \cdot \nu, \\
        \sum_{j=1}^M J_{p,j}(t) \int_\Omega \phi_j\cdot \phi_l -\sum_{i=1}^N p_i(t) \int_\Omega \varphi_i\ \divergence \phi_l + \int_\Omega p^h  g^h_{\psi} \cdot \phi_l&=- \int_\Gamma p^h\ \phi_l \cdot \nu,
    \end{align*}
which represents a nonlinear, large and sparse DAE for the approximation of the functions $\psi$, $n$, $p$, $g_\psi$, $J_n$, and $J_p$. In matrix notation it reads
    \begin{multline*}
         \begin{pmatrix}
             0 \\
            -M_L \dot{n}(t) \\
            \phantom{-} M_L \dot{p}(t) \\
             0 \\
             0 \\
             0
         \end{pmatrix}
         +
         \underbrace{\begin{pmatrix}
             & -M_L & M_L & \lambda D & &  \\
             & & & & \nu_n D &  \\
             & & & & & \nu_p D  \\
             D^\top & & & M_H & &  \\
             & D^\top & & & M_H &  \\
             & & -D^\top & & & M_H
         \end{pmatrix}}_{A_{FEM}}
         \begin{pmatrix}
            \psi(t) \\
            n(t) \\
            p(t) \\
            g_\psi(t) \\
            J_n(t) \\
            J_p(t)
         \end{pmatrix}
         \\ + \mathcal{F}(n^h, p^h, g_\psi^h )
         = b(e(t)),
    \end{multline*}
with
    \[ \mathcal{F}(n^h, p^h, g_\psi^h) := \begin{pmatrix}
            0\\
            -\int_\Omega R(n^h, p^h) \ \varphi\\
            \phantom{-}\int_\Omega R(n^h, p^h)\ \varphi\\
            0 \\
            \int_\Omega n^h  g^h_{\psi} \cdot \phi \\
            \int_\Omega p^h  g^h_{\psi} \cdot \phi
           \end{pmatrix},\quad
        b := \begin{pmatrix}
            -\int_\Omega C \ \varphi\\
            0 \\
            0 \\
            \phantom{-} \int_\Gamma \psi^h(e(t))\ \phi \cdot \nu \\
            \phantom{-}\int_\Gamma n^h\ \phi \cdot \nu \\
            -\int_\Gamma p^h\ \phi \cdot \nu
           \end{pmatrix}, \]
and
        \[ \int_\Omega R(n^h, p^h) \varphi := \begin{pmatrix}
               \int_\Omega R(n^h, p^h) \varphi_1 \\
               \vdots \\
               \int_\Omega R(n^h, p^h) \varphi_N \\
           \end{pmatrix}. \]
All other integrals in $\mathcal{F}$ and $b$ are defined analogously. The matrices $M_L \in \R^{N \times N}$ and $M_H \in \R^{M \times M}$ are mass matrices in the spaces $L_h$ and $H_{h,0,N}$, respectively, and $D \in \R^{N \times M}$. The final DAE now takes the form
\pagebreak

\begin{problem}[full model]\label{pb:full}
\begin{align}
        A_C \frac{d}{dt} q_C(A_C^\top e(t), t) + A_R g(A_R^\top e(t), t) + A_L j_L(t) + A_V j_V(t) \nonumber\\
              + A_S j_S(t) + A_I i_s(t) &= 0, \label{eq:dae:1} \\
        \frac{d}{dt} \phi_L(j_L(t), t) - A_L^\top e(t) &= 0, \label{eq:dae:2} \\
        A_V^\top e(t) - v_s(t)&= 0, \label{eq:dae:3}\\
    		j_S(t) - C_1 J_n(t) - C_2 J_p(t) - C_3 \dot{g}_\psi(t) &= 0, \label{eq:dae:4} \\
         \begin{pmatrix}
             0 \\
            -M_L \dot{n}(t) \\
            \phantom{-} M_L \dot{p}(t) \\
             0 \\
             0 \\
             0
         \end{pmatrix}
         +
         A_{FEM}
         \begin{pmatrix}
            \psi(t) \\
            n(t) \\
            p(t) \\
            g_\psi(t) \\
            J_n(t) \\
            J_p(t)
         \end{pmatrix}
         + \mathcal{F}(n^h, p^h, g_\psi^h )
         - b(e(t)) &=0, \label{eq:dae:5}
\end{align}
where \eqref{eq:dae:4} represents the discretized linear coupling condition~\eqref{eq:current}.
\end{problem}

\section{Model reduction}\label{sec:red}
We now aim to reduce the computational effort of repeated dynamical simulations by applying proper orthogonal decomposition (POD) to the drift-diffusion equations. The idea is to replace the large number of local model-independent ansatz and test functions $\{\phi_i\}, \{\varphi_j\}$ by only a few nonlocal model-dependent Ansatz functions for the respective variables.

The snapshot variant of POD introduced in \cite{sir87} works as follows. We run a simulation of the unreduced system and collect $l$ snapshots $\psi^h(t_k,\cdot)$, $n^h(t_k,\cdot)$, $p^h(t_k,\cdot)$, $g_\psi^h(t_k,\cdot)$, $J_n^h(t_k,\cdot)$, $J_p^h(t_k,\cdot)$ at time instances $t_k \in \{t_1,\ldots,t_l\} \subset [0,T]$. The optimal selection of the time instances is not considered here. We use the time instances delivered by the DAE integrator.

Since every component of the state vector $y := (\psi, n, p, g_\psi, J_n, J_p)$ has its own physical meaning we apply POD MOR to each component separately. Among other things this approach has the advantage of yielding a block-dense model and the approximation quality of each component is adapted individually.

Let $X$ denote a Hilbert space and let $y^h : [0,T] \times X \rightarrow \R^r$ with some $r \in \N$. The Galerkin formulation~\eqref{eq:galerkinansatz} yields $y^h(t,\cdot) \in X_h := \spanvec \{\phi^X_1,\dots,\phi^X_n\}$, where $\{\phi^X_j\}_{1 \leq j \leq n}$ denote $n$ linearly independent elements of $X$. The idea of POD consists in finding a basis $\{u^1, \ldots, u^m\}$ of the span of the snapshots
    \[ \spanvec\left\{ y^h(t_k,\cdot) = \sum_{i=1}^n y^{h,k}_i \phi^X_i(\cdot),\text{ with } k=1,\dots,l\right\} \]
satisfying
    \[ \{u^1,\ldots,u^s\} = \argmin_{\{v^1,\ldots,v^s\} \subset X} \quad \sum_{k=1}^l \Big\|y^h(t_k,\cdot)-\sum_{i=1}^{s}\langle y^h(t_k,\cdot),v^i(\cdot)\rangle_X v^i(\cdot)\Big\|^2_X\ , \]
for $1\leq s\leq m$, where $1\leq m \leq l$. The functions $\{u^i\}_{1 \leq i \leq s}$ are orthonormal in $X$ and can be obtained with the help of SVD as follows.

Let the matrix $Y:=(y^{h,1},\dots,y^{h,l})$ $\in \R^{n \times l}$ contain as columns the coefficient vectors of the snapshots. Furthermore, let
	$ M:= (\langle \phi^X_i,\phi^X_j\rangle_X )_{1\leq i,j \leq n} $
be the positive definite mass matrix with its Cholesky factorization $M=LL^\top$. Let $(\tilde U,\Sigma,\tilde V)$ denote the singular value decomposition of $\tilde Y := L^\top Y$, i.e. $\tilde Y=\tilde U\Sigma \tilde V^\top$ with $\tilde U\in \R^{n\times n}$, $\tilde V\in \R^{l\times l}$, and a matrix $\Sigma \in \R^{n\times l}$ containing the singular values $\sigma_1 \geq \sigma_2 \geq \ldots \geq \sigma_m > \sigma_{m+1}\geq \ldots \geq \sigma_l \geq 0$. We set $U := L^{-\top} \tilde U_{(:,\,1:s)}$. Then, the $s$-dimensional POD basis is given by
\[
\spanvec\left\{u^i(\cdot) = \sum\limits_{j=1}^n U_{ji} \phi^X_j(\cdot),\ i=1,\dots,s\right\}.
\]
The information content of $\{u^1,\ldots,u^s\}$ with respect to the scalar product $\langle \cdot, \cdot \rangle_X$ with
  \begin{equation}\label{eq:delta}
  0 \leq \Delta(s) = \sqrt{\frac{\sum_{i=s+1}^m\sigma_i^2}{\sum_{i=1}^m\sigma_i^2}} \leq 1,
  \end{equation}
is given by $1-\Delta(s)$. Here $\Delta(s)$ measures the lack of information of $\{u^1,\ldots,u^s\}$ with respect to $\spanvec\{y^h(t_1,\cdot), \ldots, y^h(t_l,\cdot)\}$.

The POD basis functions are now used as trial and test functions in the Galerkin method.

If the snapshots satisfy inhomogeneous Dirichlet boundary conditions POD is performed for
\[
 \tilde \psi(t) = \psi(t) - \psi_r(t),\quad \tilde n(t)=n(t)-n_r(t), \quad \tilde p(t) = p(t)-p_r(t),
\]
with $\psi_r$, $n_r$, $p_r$ denoting reference functions satisfying the Dirichlet boundary conditions required for $\psi$, $n$ and $p$. This guarantees that the POD basis admits homogeneous boundary conditions on the Dirichlet boundary.
In the case of the mixed finite element approach the introduction of a reference state is not necessary, since the boundary values are included more naturally through the variational formulation.

The time-snapshot POD procedure delivers Galerkin Ansatz spaces for $\psi$, $n$, $p$, $g_\psi$, $J_n$ and $J_p$. This leads to the Ansatz
\begin{equation}\label{eq:podansatz}
    \left.
    \begin{aligned}
      \psi^{POD}(t)    &= U_\psi \gamma_\psi(t), \quad
      &g_\psi^{POD}(t) &= U_{g_\psi} \gamma_{g_\psi}(t), \\
      n^{POD}(t)       &= U_n \gamma_n(t), \quad
      &J_n^{POD}(t)    &= U_{J_n} \gamma_{J_n}(t),\\
      p^{POD}(t)       &= U_p \gamma_p(t), \quad
      &J_p^{POD}(t)    &= U_{J_p} \gamma_{J_p}(t).
    \end{aligned}
    \quad \right\}
\end{equation}
The injection matrices
\begin{align*}
  &U_\psi  \in \R^{N\times s_\psi },
  &&U_n     \in \R^{N\times s_n    },
  &&U_p     \in \R^{N\times s_p    },\\
  &U_{g_\psi}  \in \R^{M\times s_{g_\psi} },
  &&U_{J_n} \in \R^{M    \times s_{J_n}},
  &&U_{J_p} \in \R^{M    \times s_{J_p}},
\end{align*}
contain the (time independent) POD basis functions, the vectors $\gamma_{(\cdot)}$ the corresponding time-variant coefficients. The numbers $s_{(\cdot)}$ are the respective number of POD basis functions included.
Assembling the POD system yields the DAE
    \begin{equation*}
         \begin{pmatrix}
             0 \\
            -\dot{\gamma}_n(t) \\
            \phantom{-} \dot{\gamma}_p(t) \\
             0 \\
             0 \\
             0
         \end{pmatrix}
         +
         A_{POD}
         \begin{pmatrix}
            \gamma_{\psi}(t) \\
            \gamma_{n}(t) \\
            \gamma_{p}(t) \\
            \gamma_{g_\psi}(t) \\
            \gamma_{J_n}(t) \\
            \gamma_{J_p}(t)
         \end{pmatrix}
          + U^\top \mathcal{F}( n^{POD}, p^{POD}, g_\psi^{POD} ) = U^\top b(e(t)),
    \end{equation*}
with
\begin{align*}
	A_{POD} &= U^\top A_{FEM} U \\
	        &= \begin{pmatrix}
             & -U_\psi^\top M_L U_n & U_\psi^\top M_L U_p & \lambda U_\psi^\top D U_{g_\psi} & &  \\
             & & & & \nu_n U_n^\top D U_{J_n}&  \\
             & & & & & \nu_p U_p^\top D U_{J_p}  \\
             U_{g_\psi}^\top D^\top U_\psi & & & I & &  \\
             & U_{J_n}^\top D^\top U_n & & & I &  \\
             & & -U_{J_p}^\top D^\top U_p & & & I
         \end{pmatrix}
\end{align*}
and $U = \diag( U_\psi, U_n, U_p, U_{g_\psi}, U_{J_n}, U_{J_p} )$. Note that we exploit the orthogonality of the POD basis functions, e.g.  $U_n^\top M_L U_n = U_p^\top M_L U_p = I_{N \times N}$ and $U_{g_\psi}^\top M_H U_{g_\psi}  = U_{J_n}^\top M_H U_{J_n} = U_{J_p}^\top M_H U_{J_p} = I_{M \times M}$. The arguments of the nonlinear functional have to be interpreted as functions in space.

All matrix-matrix multiplications are calculated in an offline-phase. The nonlinear functional $\mathcal{F}$ has to be evaluated online. A possible reduction method for the nonlinearity is called Discrete Empirical Interpolation and is discussed in~\cite{CS09}. Here, we focus on the speed-up in solving linear systems in the implicit time integrator, which are now small. The reduced model for the network now reads
\pagebreak
\begin{problem}[reduced order model]\label{pb:rom}
\begin{align}
        A_C \frac{d}{dt} q_C(A_C^\top e(t), t) + A_R g(A_R^\top e(t), t) + A_L j_L(t) + A_V j_V(t) \nonumber\\
              + A_S j_S(t) + A_I i_s(t) &= 0, \label{eq:daered:1} \\
        \frac{d}{dt} \phi_L(j_L(t), t) - A_L^\top e(t) &= 0, \label{eq:daered:2} \\
        A_V^\top e(t) - v_s(t)&= 0, \label{eq:daered:3}\\
    		j_S(t) - C_1 U_{J_n} \gamma_{J_n}(t) - C_2 U_{J_p} \gamma_{J_p}(t) - C_3 U_{g_\psi} \dot{\gamma}_{g_\psi}(t) &= 0, \label{eq:daered:4} \\
                \begin{pmatrix}
             0 \\
            -\dot{\gamma}_n(t) \\
            \phantom{-} \dot{\gamma}_p(t) \\
             0 \\
             0 \\
             0
         \end{pmatrix}
         +
         A_{POD}
         \begin{pmatrix}
            \gamma_{\psi}(t) \\
            \gamma_{n}(t) \\
            \gamma_{p}(t) \\
            \gamma_{g_\psi}(t) \\
            \gamma_{J_n}(t) \\
            \gamma_{J_p}(t)
         \end{pmatrix}
          + U^\top \mathcal{F}( n^{POD}, p^{POD}, g_\psi^{POD} ) - U^\top b(e(t)) &= 0. \label{eq:daered:5}
\end{align}
\end{problem}

\section{Residual-based sampling}\label{sec:param}
We now consider parameter dependent models. Possible parameters are physical constants of the semiconductors (e.g. length, permeability, doping) and parameters of the network elements (e.g. frequency of sinusoidal voltage sources, value of resistances). We do not distinguish between inputs and parameters of the model.

Let there be $r \in \N$ parameters and let the space of considered parameters be given as a bounded set $\calP \subset \R^r$. We construct the reduced model based on snapshots from a simulation at a reference parameter $\omega_1 \in \calP$. One expects that the reduced model approximates the unreduced model well in a small neighborhood of $\omega_1$, but one cannot expect that the reduced model is valid over the complete parameter set $\calP$. In order to create a suitable reduced order model we consider additional snapshots which are obtained from simulations at parameters $\omega_2,\omega_3,\ldots \in \calP$. The iterative selection of $\omega_{k+1}$ at a step $k$ is called parameter sampling. Let $P_k$ denote the set of selected reference parameters, $P_k := \{ \omega_1, \omega_2, \ldots, \omega_k\} \subset \calP$.

We neglect the discretization error of the finite element method and its influence on the coupled network and define the error of the reduced model as
	\begin{equation}\label{eq:deferror}
			\mathcal{E}(\omega; P) := z^h(\omega) - z^{POD}(\omega; P),
	\end{equation}
where $z^h(\omega) := ( e^h(\omega), j_V^h(\omega), j_L^h(\omega), y^h(\omega) )^\top$ is the solution of problem~\ref{pb:full} at the parameter $\omega$ with discretized semiconductor variables $y^h := (\psi^h, n^h, p^h, g_\psi^h, J_n^h, J_p^h)^\top$. $z^{POD}(\omega; P)$ denotes the solution of the coupled system in problem~\ref{pb:rom} with reduced semiconductors, where the reduced model is created based on simulations at the reference parameters $P \subset \calP$. The error is considered in the space $X$ with norm
	\begin{align*}
			\| z \|_X :=\Big\| \Big( & \| e \|_2, \| j_V \|_2, \| j_L \|_2, \\
			& \| \psi \|_{L^2([0,T], L^2(\Omega))}, \| n \|_{L^2([0,T], L^2(\Omega))}, \| p \|_{L^2([0,T], L^2(\Omega))}, \\
			& \| g_\psi \|_{L^2([0,T], H_{0,N}(\divergence, \Omega))},\\
            & \| J_n \|_{L^2([0,T], H_{0,N}(\divergence, \Omega))}, \| J_p \|_{L^2([0,T], H_{0,N}(\divergence, \Omega))} \Big) \Big\|.
	\end{align*}
Obvious extensions apply when there is more than one semiconductor present.

Furthermore we define the residual $\mathcal{R}$ by evaluation of the unreduced model~\eqref{eq:dae:1}-\eqref{eq:dae:5} at the solution of the reduced model $z^{POD}(\omega;P)$, i.e.
	\begin{multline}\label{eq:defresidual}
			\mathcal{R}(z^{POD}(\omega;P)) :=
			   \begin{pmatrix}
             0 \\
            -M_L \dot{n}^{POD}(t) \\
            \phantom{-} M_L \dot{p}^{POD}(t) \\
             0 \\
             0 \\
             0
         \end{pmatrix}
         +
         A_{FEM}
         \begin{pmatrix}
            \psi^{POD}(t) \\
            n^{POD}(t) \\
            p^{POD}(t) \\
            g_\psi^{POD}(t) \\
            J_n^{POD}(t) \\
            J_p^{POD}(t)
         \end{pmatrix} \\
         + \mathcal{F}(n^{POD}, p^{POD}, g_\psi^{POD} )
         - b(e^{POD}(t)).
	\end{multline}
Note that the residual of equations~\eqref{eq:dae:1}-\eqref{eq:dae:4} vanishes.

We note that the same definitions are used in~\cite{HO09} for linear descriptor systems. In~\cite{HO09} an error estimate is obtained by deriving a linear ODE for the error and exploiting explicit solution formulas. Here we have a nonlinear DAE and at the present state we are not able to provide an upper bound for the error $\|\mathcal{E}(\omega;P)\|_X$ which would yield a rigorous sampling method using for example the Greedy algorithm of~\cite{PR07}.

We propose to consider the residual as an estimate for the error. The evaluation of the residual is cheap since it only requires the solution of the reduced system and its evaluation in the unreduced DAE. It is therefore possible to evaluate the residual at a large set of test parameters $P_{test} \subset \calP$. Similar to the Greedy algorithm of~\cite{PR07}, we add to the set of reference parameters the parameter where the residual becomes maximal.

The magnitude of the components in error and residual may be large and a proper scaling should be applied. For the error we consider the component-wise relative error, i.e.
\[ \frac{\| \psi^h(\omega) - \psi^{POD}(\omega; P) \|_{L^2([0,T], L^2(\Omega))}}{\| \psi^h(\omega) \|_{L^2([0,T], L^2(\Omega))}},\
   \frac{\| n^h(\omega) - n^{POD}(\omega; P) \|_{L^2([0,T], L^2(\Omega))}}{\| n^h(\omega) \|_{L^2([0,T], L^2(\Omega))}},\ \ldots,
\]
and the residual is scaled by a block-diagonal matrix containing the weights
\begin{multline*}
    D(\omega) \mathcal{R}(z^{POD}(\omega;P)) =\\
    \left(\begin{array}{@{}c@{}c@{}c@{}c@{}c@{}c@{}}
        d_\psi(\omega) I \\
        & d_n(\omega) I \\
        & & d_p(\omega) I\\
        & & & d_{g_\psi}(\omega) I\\
        & & & & d_{J_n}(\omega) I \\
        & & & & & d_{J_p}(\omega) I
    \end{array}\right)
 \mathcal{R}(z^{POD}(\omega;P)).
\end{multline*}
The weights $d_{(\cdot)}(\omega)>0$ may be parameter-dependent. These weights are chosen in a way that the norm of the residual and the relative error are component-wise equal at the reference frequencies $\omega_k$ where we know $z^h(\omega_k)$ from simulation of the unreduced model, i.e.
\begin{equation}\label{eq:weight:interp}
    d_\psi(\omega_k) := \frac{\| \psi^h(\omega_k) - \psi^{POD}(\omega_k; P) \|_{L^2([0,T], L^2(\Omega))}}{\| \psi^h(\omega_k) \|_{L^2([0,T], L^2(\Omega))} \cdot \|\mathcal{R}_1(z^{POD}(\omega_k;P))\|_{L^2([0,T], L^2(\Omega))}},
\end{equation}
and similarly for the other components. If $\|\mathcal{R}_1(z^{POD}(\omega_k;P))\|_{L^2([0,T], L^2(\Omega))} = 0$ we chose $d_\psi(\omega_k) := 1$.

In 1D-parameter sampling with $\calP := [ \underline{p},\ \overline{p} ]$, this can be achieved by interpolating the weights $d_{(\cdot)}(\omega_1)$, $\ldots$, $d_{(\cdot)}(\omega_k)$ piecewise linearly. Extrapolation is done by nearest-neighbour interpolation to ensure the positivity of the weights.

We summarize our ideas in the following sampling algorithm:
\begin{algorithm}[Sampling]\label{algorithm:sampling}\
\begin{enumerate}
\item Select $\omega_1 \in \calP$, $P_{test} \subset \calP$, $tol > 0$, and set $k:=1$, $P_1 := \{ \omega_1 \}$.
\item Simulate the unreduced model at $\omega_1$ and calculate the reduced model with POD basis functions $U_1$.
\item Calculate weight functions $d_{(\cdot)}(\omega) > 0$ according to \eqref{eq:weight:interp} for all $\omega_k \in P_k$.
\item Calculate the scaled residual $\|D(\omega) \mathcal{R}(z^{POD}(\omega, P_k))\|$ for all $\omega \in P_{test}$.
\item Check termination conditions, e.g.
\begin{itemize}
 \item $\max_{\omega \in P_{test}} \|D(\omega) \mathcal{R}(z^{POD}(\omega, P_k))\| < tol$,
 \item no progress in weighted residual.
\end{itemize}
\item Calculate $\omega_{k+1} := \argmax_{\omega \in P_{test}} \|D(\omega)\mathcal{R}(z^{POD}(\omega, P_k))\|$.
\item Simulate the unreduced model at $\omega_{k+1}$ and create a new reduced model with POD basis $U_{k+1}$ using also the already available information at $\omega_1$, $\ldots$, $\omega_k$.
\item Set $P_{k+1} := P_k \cup \{ \omega_{k+1} \}$, $k := k + 1$ and goto 3.
\end{enumerate}
\end{algorithm}
The step 7 in Algorithm~\ref{algorithm:sampling} can be executed in different ways. If offline time and offline memory requirements are not critical one may combine snapshots from all simulations of the full model and redo the model order reduction on the large snapshot ensemble. Otherwise we can create a new reduced model at reference frequency $\omega_{k+1}$ with POD-basis $\bar{U}$ and then performing an additional POD step on $( U_k, \bar{U})$.

\section{Numerical investigation}\label{sec:numerics}
The discussed finite element method is implemented for spatial one-dimensional (1D) problems in MATLAB. The resulting high dimensional DAE is integrated using the DASSL software package~\cite{Pet93}. We assume that the differentiation index of the network is $1$. Otherwise one should switch to alternative integrators. In order to solve the Newton systems efficiently which arise from the BDF method, the variables of the sparse system are reordered with respect to minimal bandwidth.

\begin{figure}
\includegraphics[width=.322\textwidth]{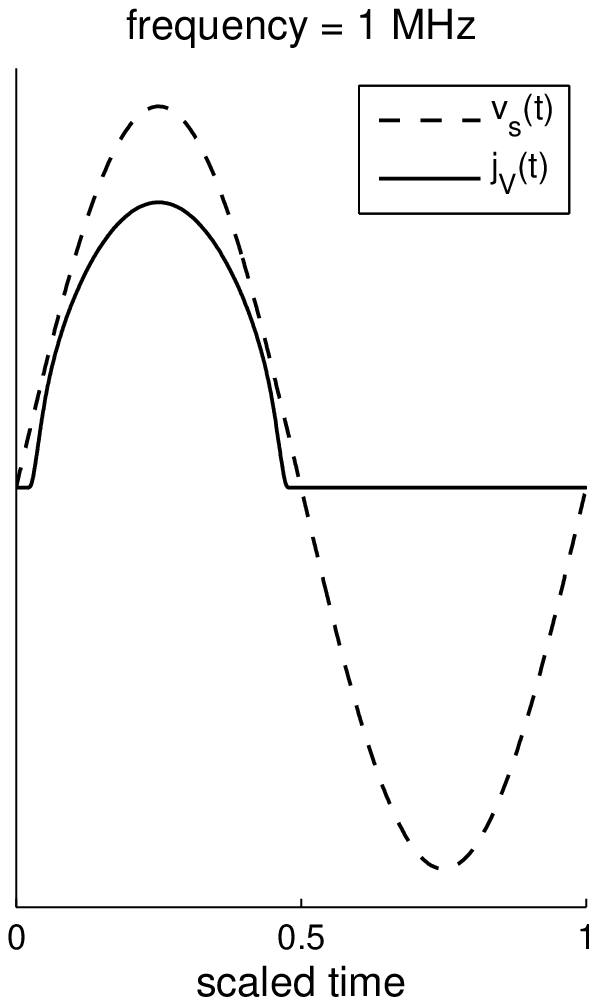}
\includegraphics[width=.322\textwidth]{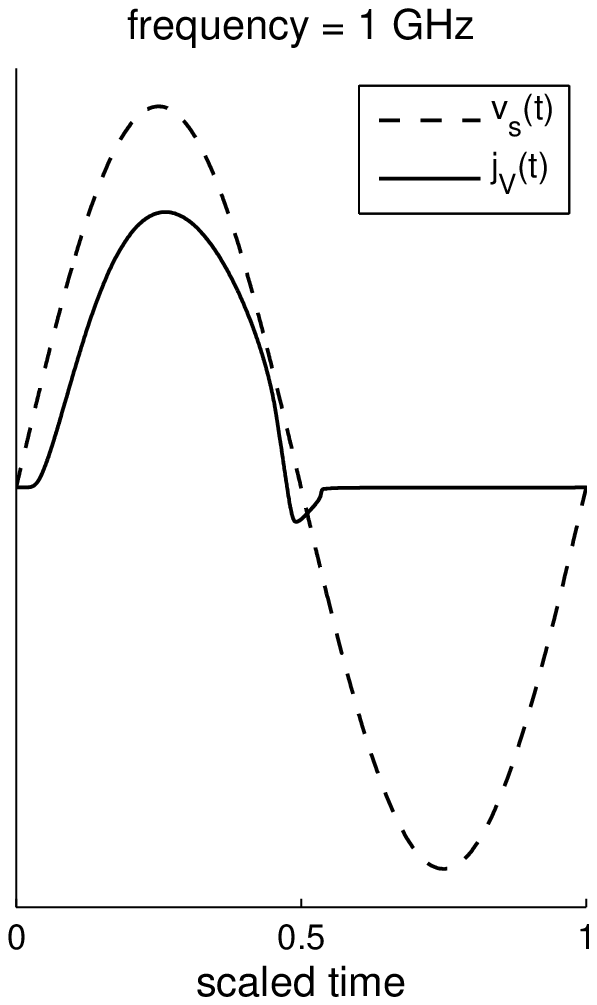}
\includegraphics[width=.322\textwidth]{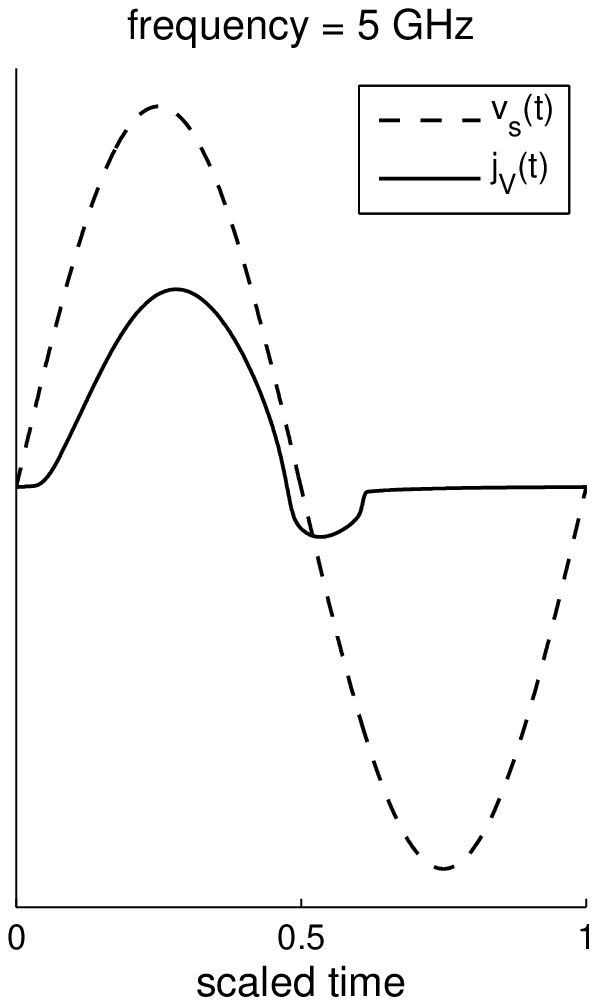}

\caption{Current $j_V$ through the basic network for input frequencies 1 MHz, 1 GHz and 5 GHz. The capacitive effect is clearly demonstrated.}
\label{fig:basicres}
\end{figure}

A basic test circuit with a single diode is depicted in Figure \ref{fig:basic}, where the model parameters are presented in Table \ref{tab:diodeparam}. The input $v_s(t)$ is chosen to be sinusoidal with amplitude $5\ [V]$. The numerical results in Figure \ref{fig:basicres} show the capacitive effect of the diode for high input frequencies. Similar results are obtained in \cite{Sot06} using the simulator MECS. In the sequel the frequency will be considered as a model parameter.

\begin{table}
\caption{Diode model parameters.}
\label{tab:diodeparam}
\begin{tabular}{p{1.9cm}p{1.0cm}p{2.40cm}p{2.0cm}p{2.1cm}p{1.6cm}}
\hline\noalign{\smallskip}
Parameter & Value & & Parameter & Value & \\
\noalign{\smallskip}\svhline\noalign{\smallskip}
$ L$ (length) & $ 10^{-4} $ & [$cm$]           & $ \varepsilon   $ & $ 1.03545 \cdot 10^{-12} $ & [$F/cm$]    \\
$ U_T   $ & $ 0.0259  $ & [$V$]            & $ \eta          $ & $ 1.4\cdot 10^{10}       $ & [$1/cm^3$]  \\
$ \mu_n $ & $ 1350    $ & [$cm^2/(V sec)$] & $ \tau_n        $ & $ 330\cdot 10^{-9}       $ & [$sec$]     \\
$ \mu_p $ & $ 480     $ & [$cm^2/(V sec)$] & $ \tau_p        $ & $ 33 \cdot 10^{-9}       $ & [$sec$]     \\
$ a$ (contact area) & $ 10^{-5} $ & [$cm^2$]         & $ C(x),\ x<L/2   $ & $ -9.94 \cdot 10^{15}    $ & [$1/cm^3$]  \\
          &             &                & $ C(x),\ x\geq L/2 $ & $ 4.06 \cdot 10^{18}     $ & [$1/cm^3$]  \\
\noalign{\smallskip}\hline\noalign{\smallskip}
\end{tabular}
\end{table}

\begin{figure}[t]
\sidecaption
\includegraphics[width=0.55\textwidth]{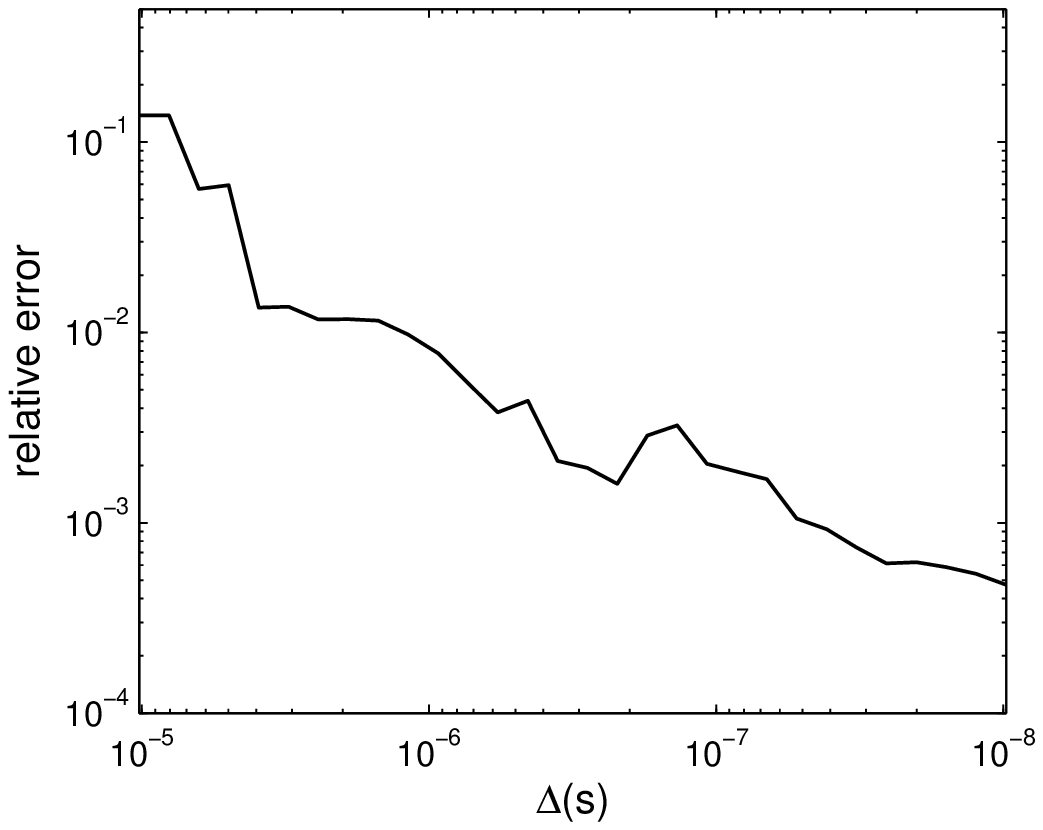}
\caption{Relative error between reduced and unreduced problem at the fixed frequency $10^{10}\ [Hz]$.}
\label{fig:err_over_delta}
\end{figure}

We first validate the reduced model at a fixed reference frequency of $10^{10}\ [Hz]$. Figure \ref{fig:err_over_delta} shows the development of the relative error between the reduced and the unreduced numerical solutions, plotted over the lack of information $\Delta$, see \eqref{eq:delta}. The number of POD basis functions for each variable is chosen such that the indicated approximation quality is reached, i.e.
$
 \Delta:=\Delta_\psi\simeq\Delta_n\simeq\Delta_p\simeq\Delta_{g_\psi}\simeq\Delta_{J_n}\simeq\Delta_{J_p}\,.
$
Since we compute all POD basis functions anyway, this procedure does not involve any additional costs.

\begin{figure}[t]
\sidecaption
\includegraphics[width=0.55\textwidth]{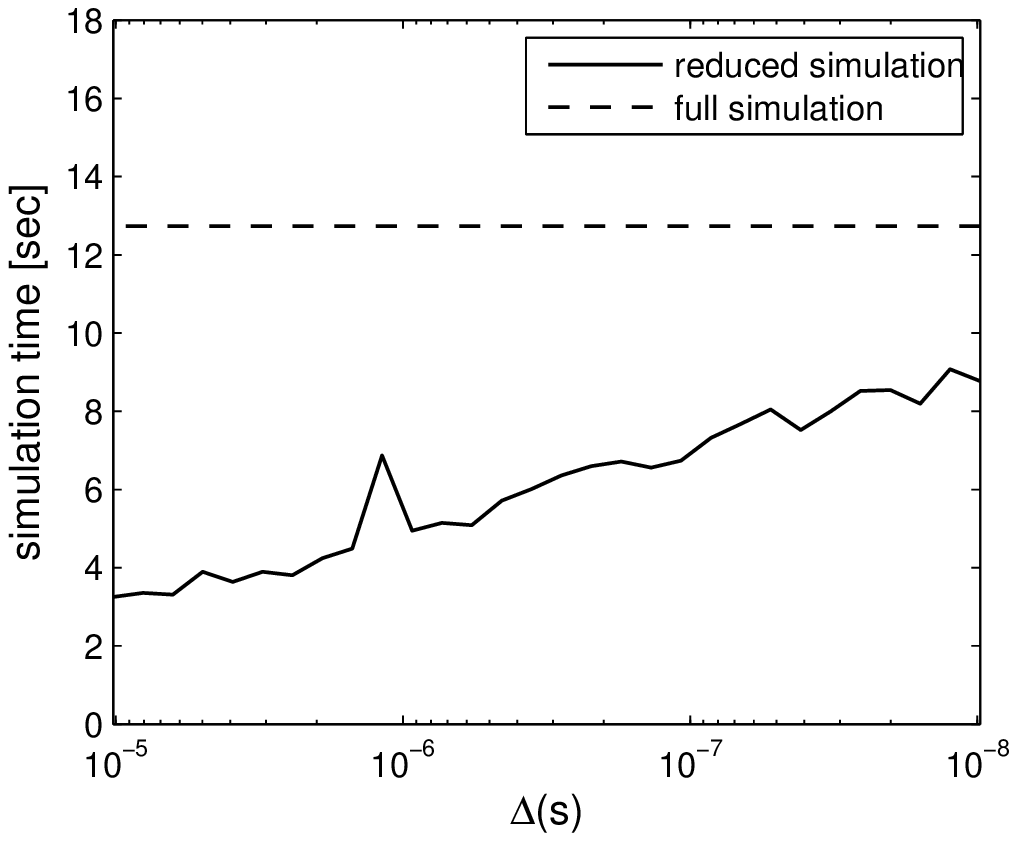}
\caption{Time consumption for simulation runs of Figure~\ref{fig:err_over_delta}. The dashed line indicates the time consumption for the simulation of the original full system.}
\label{fig:time_over_delta}
\end{figure}

In Figure \ref{fig:time_over_delta} the simulation times are plotted versus the lack of information $\Delta$.

\begin{figure}[t] \sidecaption
\includegraphics[width=0.55\textwidth]{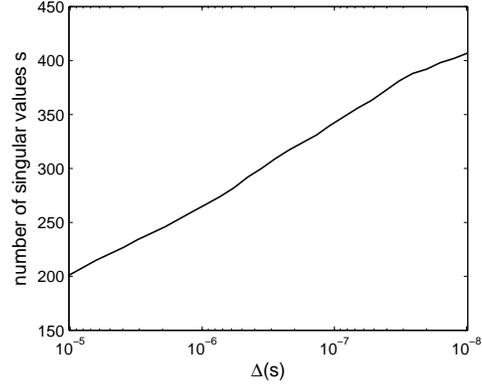}
\caption{The number of required singular values grows only logarithmically with the requested information content.}
\label{fig:singvals}
\end{figure}

Figure~\ref{fig:singvals} shows the total number of singular vectors $s$ required to guarantee a certain information content $1-\Delta$. It can be seen that with the number of singular vectors included increasing only linearly, the lack of information tends to zero exponentially.

\begin{figure}[h]
\sidecaption
\includegraphics[width=0.55\textwidth]{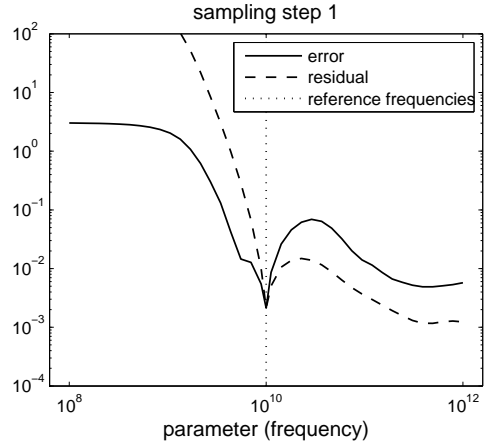}
\caption{Relative reduction error (solid line) and weighted residual (dashed line) plotted over the frequency parameter space. The reduced model is created based on simulations at the reference frequency $\omega_1 := 10^{10}\ [Hz]$. The reference frequencies are marked by vertical dotted lines.}
\label{fig:sampling:1}
\end{figure}

We now apply Algorithm~\ref{algorithm:sampling} to provide a reduced order model of the basic circuit and we choose the frequency of the input voltage $v_s$ as model parameter. As parameter space we chose the interval $\calP := [10^8,\ 10^{12}]\ [Hz]$. We start the investigation with a reduced model which is created from the simulation of the full model at the reference frequency $\omega_1 := 10^{10}\ [Hz]$. The number of POD basis functions $s$ is chosen such that the lack of information $\Delta(s)$ is approximately $10^{-7}$. The relative error and the weighted residual are plotted in Figure~\ref{fig:sampling:1}. We observe that the weighted residual is a rough estimate for the relative approximation error. Using Algorithm~\ref{algorithm:sampling} the next additional reference frequency is $\omega_2 := 10^{8}\ [Hz]$ since it maximizes the weighted residual. The second reduced model is constructed on the same lack of information $\Delta := 10^{-7}$. Here we note that in practical applications, the error is not known over the whole parameter space.

\begin{figure}[t]
\sidecaption
\includegraphics[width=0.55\textwidth]{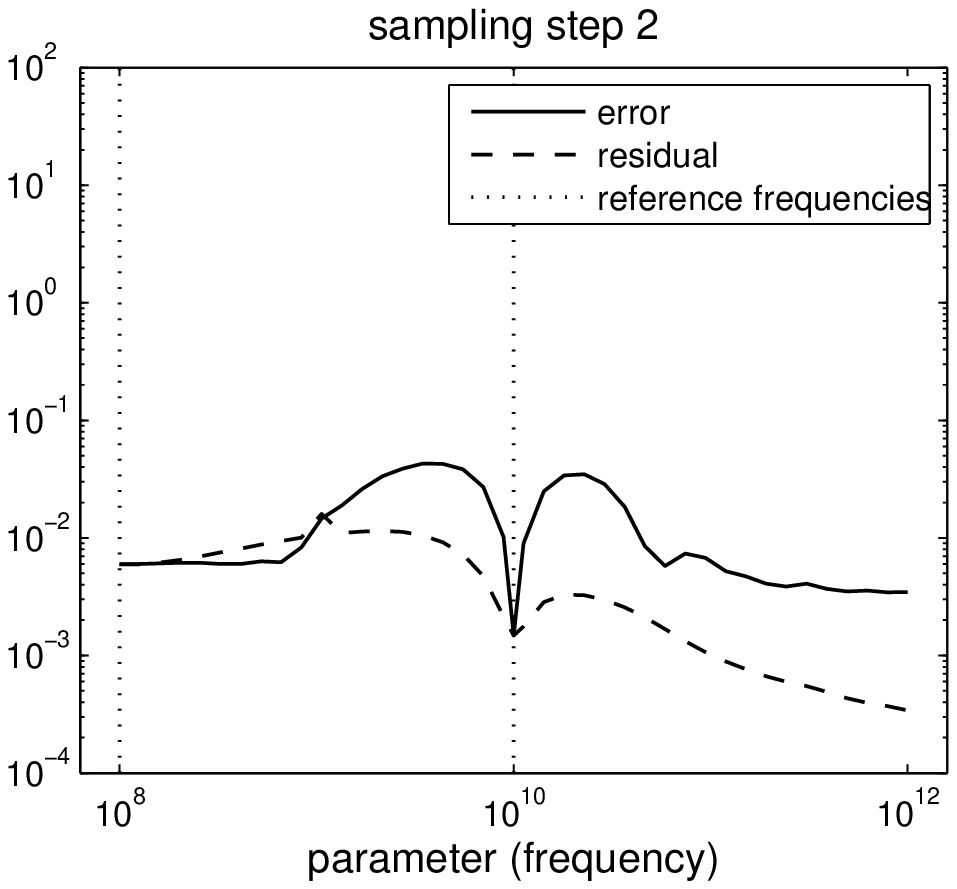}
\caption{Relative reduction error (solid line) and weighted residual (dashed line) plotted over the frequency parameter space. The reduced model is created based on simulations at the reference frequencies $\omega_1 := 10^{10}\ [Hz]$ and $\omega_2 := 10^{8}\ [Hz]$. The reference frequencies are marked by vertical dotted lines.}
\label{fig:sampling:2}
\end{figure}

\begin{figure}[t]
\sidecaption
\includegraphics[width=0.55\textwidth]{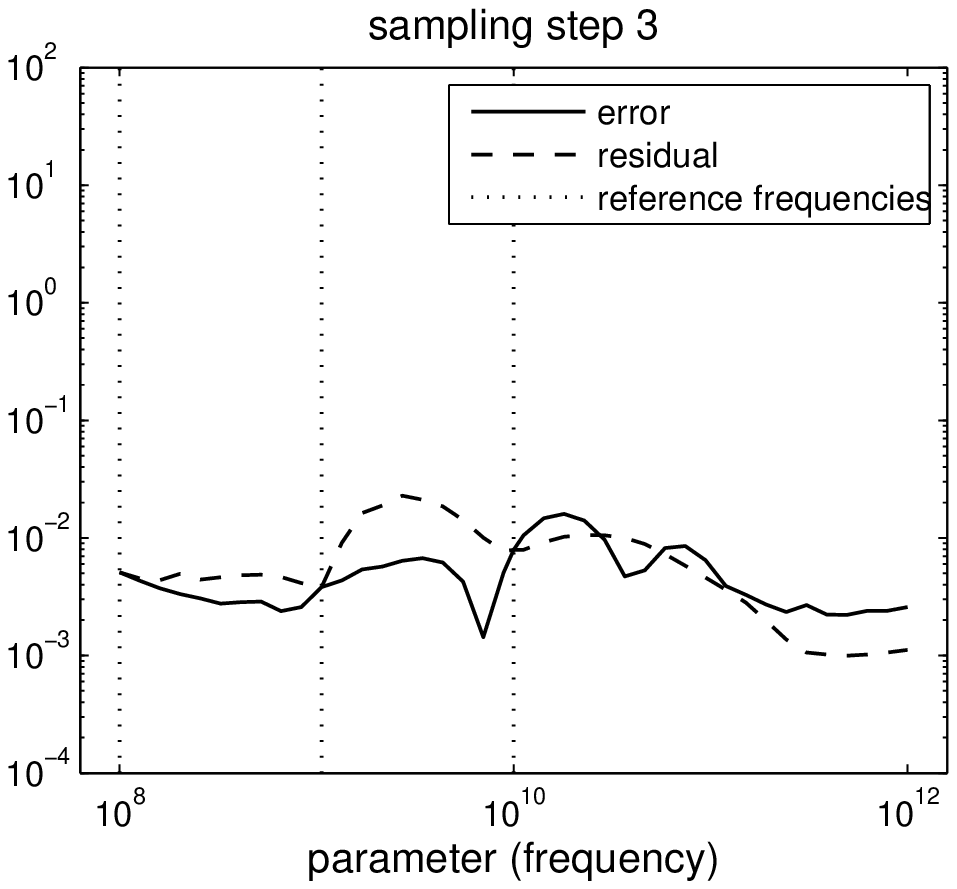}
\caption{Relative reduction error (solid line) and weighted residual (dashed line) plotted over the frequency parameter space. The reduced model is created based on simulations at the reference frequency $\omega_1 := 10^{10}\ [Hz]$, $\omega_2 := 10^{8}\ [Hz]$, and $\omega_3 := 1.0608 \cdot 10^{9}\ [Hz]$. The reference frequencies are marked by vertical dotted lines.}
\label{fig:sampling:3}
\end{figure}

The next two iterations of the sampling algorithm are depicted in Figures~\ref{fig:sampling:2} and~\ref{fig:sampling:3}. Based on the residual in step 2, one selects $\omega_3 := 1.0608 \cdot 10^{9}\ [Hz]$ as the next reference frequency. Since no further progress of the weighted residual is achieved in step 3, the algorithm terminates. The maximal errors and residuals are given in Table~\ref{table:sampling}.

\begin{table}[b]
\caption{Progress of refinement method.}
\label{table:sampling}
\begin{tabular}{p{1.0cm}p{3.5cm}p{3.0cm}p{3.0cm}}
\hline\noalign{\smallskip}
step $k$ & reference parameters $P_k$ & max. scaled residual & max. relative error \\
         &                            & (at frequency)       & (at frequency) \\
\noalign{\smallskip}\svhline\noalign{\smallskip}
1 & $ \{ 1.0000 \cdot 10^{10}\}          $ & $9.9864 \cdot 10^{2}$   & $3.2189 \cdot 10^{0}$ \\
  &                                        & \textbf{$(1.0000 \cdot 10^{8})$} & $(1.0000 \cdot 10^{8})$ \\
\noalign{\smallskip}\noalign{\smallskip}
2 & $ \{1.0000 \cdot 10^{8},             $ & $1.5982 \cdot 10^{-2}$  & $4.3567 \cdot 10^{-2}$\\
  & $ \phantom{\{}1.0000 \cdot 10^{10}\} $ & $(1.0608 \cdot 10^{9})$ & $(3.4551 \cdot 10^{9})$ \\
\noalign{\smallskip}\noalign{\smallskip}
3 & $ \{1.0000 \cdot 10^{8},             $ & $2.2829 \cdot 10^{-2}$  & $1.6225 \cdot 10^{-2}$\\
  & $ \phantom{\{}1.0608 \cdot 10^{9},   $ & $(2.7283 \cdot 10^{9})$ & $(1.8047 \cdot 10^{10})$ \\
  & $ \phantom{\{}1.0000 \cdot 10^{10}\} $ &                         & \\
\noalign{\smallskip}\hline\noalign{\smallskip}
\end{tabular}
\end{table}

\begin{figure}
\sidecaption
\includegraphics[width=0.55\textwidth]{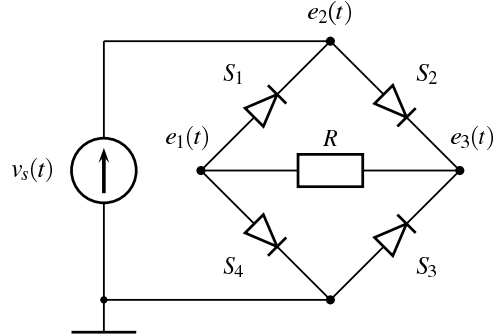}
\caption{Rectifier network.}
\label{fig:rectifier}
\end{figure}

\begin{figure}
\sidecaption
\includegraphics[width=0.55\textwidth]{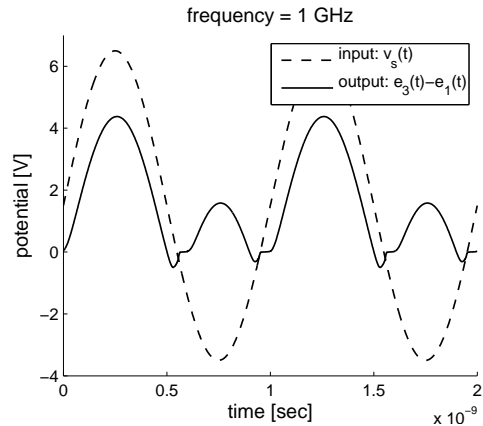}
\caption{Simulation results for the rectifier network. The input $v_s$ is sinusoidal with frequency $1\ [GHz]$ and offset $+1.5\ [V]$.}
\label{fig:rectifier2}
\end{figure}

Finally we note that the presented reduction method accounts for the position of the semiconductors in a given network in that it provides reduced order models which for identical semiconductors may be different depending on the location of the semiconductors in the network. The POD basis functions of two identical semiconductors may be different due to their different operating states. To demonstrate this fact, we consider the rectifier network in Figure~\ref{fig:rectifier}. Simulation results are plotted in Figure~\ref{fig:rectifier2}.
The distance between the spaces $U^1$ and $U^2$ which are spanned, e.g., by the POD-functions $U^1_\psi$ of the diode $S_1$ and $U^2_\psi$ of the diode $S_2$ respectively, is measured by
  \[ d(U^1, U^2) := \max_{\substack{u \in U^1\\\|u\|_2 = 1}}\,
                    \min_{\substack{v \in U^2\\\|v\|_2 = 1}}\,
                    \| u-v \|_2. \]
Exploiting the orthonormality of the bases $U^1_\psi$ and $U^2_\psi$ and using a Lagrange framework, we find
  \[ d(U^1, U^2) = \sqrt{ 2 - 2 \sqrt{ \lambda } }, \]
where $\lambda$ is the smallest eigenvalue of the positive definite matrix $S S^\top$ with $S_{ij} = \langle u^1_{\psi,i}, u^2_{\psi,j} \rangle_2$. The distances for the rectifier network are given in Table~\ref{tab:distpod}. While the reduced model for the diodes $S_1$ and $S_3$ are almost equal, the models for the diodes $S_1$ and $S_2$ are significantly different. Similar results are obtained for the reduction of $n$, $p$, etc.

\begin{table}
\caption{Distances between reduced models in the rectifier network.}
\label{tab:distpod}
\begin{tabular}{p{2.5cm}p{2.5cm}p{2.5cm}}
\hline\noalign{\smallskip}
$\Delta$  & $d(U^1,U^2)$ & $d(U^1,U^3)$ \\
\noalign{\smallskip}\svhline\noalign{\smallskip}
$10^{-4}$ & $0.61288$ & $5.373 \cdot 10^{-8}$ \\
$10^{-5}$ & $0.50766$ & $4.712 \cdot 10^{-8}$ \\
$10^{-6}$ & $0.45492$ & $2.767 \cdot 10^{-7}$ \\
$10^{-7}$ & $0.54834$ & $1.211 \cdot 10^{-6}$ \\
\noalign{\smallskip}\hline\noalign{\smallskip}
\end{tabular}
\end{table}

\begin{acknowledgement}
The work reported in this paper was supported by the German Federal Ministry of Education and Research (BMBF), grant no. 03HIPAE5. Responsibility for the contents of this publication rests with the authors.
\end{acknowledgement}

\bibliographystyle{spmpsci}

\end{document}